\RequirePackage{pgf}

\documentclass[bj,authoryear]{imsart}

\usepackage{tikz}

\usepackage{float}
\usetikzlibrary{patterns}
\usepackage{adjustbox}
\usepackage{bm}

\usepackage[textsize=scriptsize,textwidth=2.5cm]{todonotes}

\usepackage[latin1]{inputenc}

\input{macromine01.tex}

\definecolor{indigo}{RGB}{75,0,130}

\newcommand{\leqlex}{\leq_{\mathrm{lex}}}
\newcommand{\lelex}{<_{\mathrm{lex}}}
\newcommand{\leqcw}{\leq_{\mathrm{cw}}}

\newcommand{\WD}{\bm{\mathrm{WD}}}

\usepackage{hyperref}
\usepackage{cleveref}
\hypersetup{colorlinks,
	linkcolor=blue,
	citecolor=blue,
	urlcolor=magenta,
	linktocpage,
	plainpages=false}

\startlocaldefs

\newtheorem{theorem}{Theorem}
\newtheorem{lemma}[theorem]{Lemma}
\newtheorem{proposition}[theorem]{Proposition}
\newtheorem{corollary}[theorem]{Corollary}

\newcommand{\SG}{\mathrm{SG}} 

\renewenvironment{proof}
{\par\noindent{\bfseries\upshape Proof\ }}
{\hfill\ensuremath{\square}}

\theoremstyle{remark}

\newtheorem{remark}{Remark}
\newtheorem{assmpt}{Assumption}

\endlocaldefs

\begin{document}

\begin{frontmatter}
\title{Moment inequalities for sums of weakly dependent random fields}
\runtitle{Weakly dependent random fields}

\begin{aug}
\author[A]{\inits{F.}\fnms{Gilles}~\snm{Blanchard}\ead[label=e1]{gilles.blanchard@universite-paris-saclay.fr}}
\author[B]{\inits{S.}\fnms{Alexandra}~\snm{Carpentier}\ead[label=e2]{carpentier@uni-potsdam.de}}
\author[C]{\inits{T.}\fnms{Oleksandr}~\snm{Zadorozhnyi}\ead[label=e3]{oleksandr.zadorozhnyi@tum.de}}
\address[A]{ Institut de Math\'ematiques d'Orsay ,
Universit\'e Paris-Saclay, Paris, France\printead[presep={,\ }]{e1}}

\address[B]{Institut f\"ur Mathematik, Universit\"at Potsdam, Germany \printead[presep={,\ }]{e2}}

\address[C]{Lerhstuhl f\"ur Mathematische Statistik,
	TUM School of Computation, Information and Technology
	Technical University of Munich, Germany \printead[presep={,\ }]{e3}}
\end{aug}

\begin{abstract}
We derive both Azuma-Hoeffding and Burkholder-type inequalities for partial sums over a rectangular
grid of dimension $d$ of a random field satisfying a weak dependency assumption of projective type:
the difference between the expectation of an element of the random field and its conditional expectation
given the rest of the field at a distance more than $\delta$ is bounded, in $L^p$
distance, by a known decreasing function of $\delta$. The analysis is based on the combination of a multi-scale approximation of random sums by martingale difference sequences, and
of a careful decomposition of the domain. The obtained results 
extend previously known bounds 
under comparable hypotheses, and do not use the assumption of commuting  filtrations.
\end{abstract}


\end{frontmatter}

\section{Introduction}

Let $(X_t)_{t \in \mbz^d}$ be a real-valued integrable random field over the probability space $\paren{\Omega,\cF,\mbp}$,
$\cR$ be a hyperrectangular domain of $\mbz^d$ and define $S_\cR:= \sum_{t \in \cR}(X_t-\e{X_t})$ the deviation of the process sum over $\cR$ from its mean.
In this work we are interested in upper bounding the moments of $S_\cR$, under a weak dependency
assumption of projective type, expressed as follows, for some $p \in [2,\infty]$:
\begin{assmpt}[$\WD(p)$]
	\label{def:weak_spat_dep_new}
	For $t=(t_1,\ldots,t_d)\in \mbz^d$ and $k \in \mbn_{>0}$, define the $\sigma$-algebra
	\[ \cM_{t,k} \bydef \fS \{ X_{u}: u=(u_1,\ldots,u_d) \in \mbz^d,\;\,\,
	\sup_{i}(t_i - u_i)\geq k
	\},\] 
	it holds
	\begin{equation}
		\label{eq:weakdepdistnew}
		\forall t \in \mbz^d, \forall r>0 :  \norm{ \e{X_t | \cM_{t,r} } - \e{X_t}}_{p} \leq M_p \varphi_p(r),
	\end{equation}
	where $\varphi_{p}(\cdot)$ is a non-increasing function such that $\varphi_{p}(0)=1$, and $M_p$ is a real constant.
\end{assmpt}
    {
    The sigma-field $\cM_{t,k}$ corresponds to the sigma-field $\mbf{F}_{t-k}$ introduced in \cite{Basu:79} for martingale random fields.}

\begin{figure}[H]
\hfill
{
	\begin{tikzpicture}[x=1cm,y=1cm,scale=0.7]
		\draw[-stealth] (-0.5,0)--(6.5,0) node[right]{\huge $t_1$}; 
		\draw[-stealth] (0,-0.5)--(0,6.5) node[above]{\huge $t_2$}; 
		\fill[pattern=north west lines, pattern color=blue,fill opacity=1] (0,2) rectangle (2,6);
		\fill[pattern=north west lines, pattern color=blue, fill opacity=1] (0,0) rectangle 	(6,2);
		\fill (3,3)  circle[radius=3pt] node[above right,] { \huge $X_t$};
        \draw[>=latex, <->, thin] (3,2.9)-- node[right] {\large $k$} (3,2);
        \draw[thin,>=latex,<->] (2,3)-- node[above] {\large $k$} (2.9,3);
	\end{tikzpicture}
 }  
 \hfill
{	
        \begin{tikzpicture}[x=1cm,y=1cm,scale=0.7]
		\thispagestyle{empty}
		\draw[-stealth] (-0.5,0)--(6.5,0) node[right]{\huge $t_1$}; 
		\draw[-stealth] (0,-0.5)--(0,6.5) node[above]{\huge $t_2$}; 
		\fill[pattern=north west lines, pattern color=blue] (0,2) rectangle (2,6);
		\fill[pattern=north west lines, pattern color=blue] (2,4) rectangle (6,6);
		\fill[pattern=north west lines, pattern color=blue] (4,4) rectangle (6,2);
		\fill [pattern=north west lines, pattern color=blue] (0,0) rectangle 	(6,2);
		\fill (3,3)  circle[radius=3pt] node[above right,] { \huge $X_t$};
        \draw[thin,>=latex,<->] (3,2.9)-- node[right] {\large $k$} (3,2);
        \draw[thin,>=latex,<->] (2,3)-- node[above] {\large $k$} (2.9,3);
	\end{tikzpicture}
 }  
 \hfill~
\label{fig:sigma_fields}
\caption{Illustration of the $\sigma$-algebras $\mathcal M_{t,k}$ (left picture) and $\mathcal M'_{t,k}$ (right picture). The blue hashed surface represents the entries generating the $\sigma$-algebra.}
\end{figure}
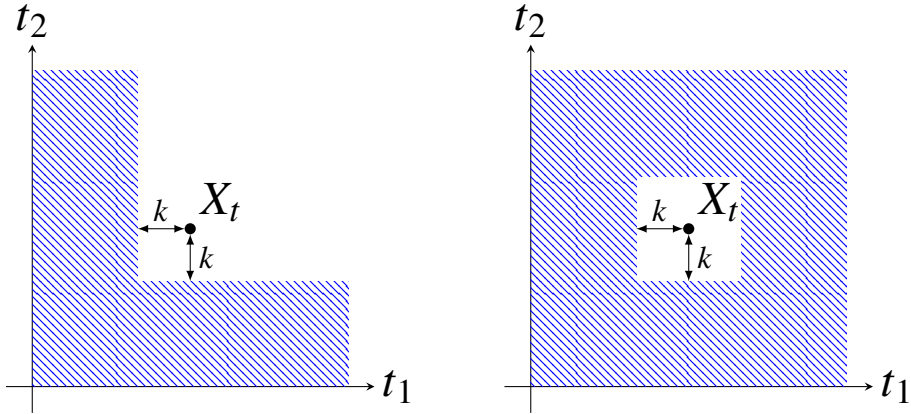

\begin{remark}
    \label{rmk:var_charactr}
    {The LHS of Equation \eqref{eq:weakdepdistnew} can be also seen as a type of generalized "variational" dependence measure between $X_{t}$ and $\cM_{t,k}$, and can be related to families of
    recently introduced covariance-based dependence coefficients (see \cite{Dedecker:05},\cite{Deschamps:06}). Namely, let $Y$ be a real random variable over $\paren{\Omega,\cF,\mbp}$ and let $\cA \subset \cF$ be some $\sigma-$field. Then for $p\geq 1$ it holds: 
    \begin{align*}
        \norm{\ee{}{Y|\cA} - \ee{}{Y}}_{p} = \sup\set{ \cov{X,Y}: \norm[1]{X}_{q} = 1, X \text{ is } \cA-\text{ measurable}}, 
    \end{align*}
    where 
    $p^{-1}+q^{-1} =1$. To see this,
    notice that $\cov{X,Y} = \e{XU}$ 
    where $U:=\e{Y|\cA}-\e{Y}$ and apply
    H\"older's inequality with extremal equality characterization (see also Lemma 1.1.2 in \cite{Deschamps:06}).
    }
\end{remark}
\begin{remark}
	 A sufficient condition for assumption $\WD(p)$ is to have~\eqref{eq:weakdepdistnew}
but with $\cM_{t,k}$ replaced by 
\begin{equation} \label{eq:mprime}
	\cM_{t,k}' \bydef \fS \{ X_{u}: u \in \mbz^d, \norm{u-t}_{\infty} \geq k \},
\end{equation}
where $\norm{\cdot}_{\infty}$ is 
the standard $\ell_{\infty}$ norm on $\mbr^d$. 
 The
$\sigma$-algebra $\cM_{t,k}'$ may actually be more intuitive or natural (see Figure~\ref{fig:sigma_fields} for an illustration), except in dimension $d=1$ where $\cM_{t,k}$ coincides with the ``past at distance $k$'' which is a natural notion for sequential processes.
\end{remark}

We will say that the random field is polynomially weakly dependent if Equation~\eqref{eq:weakdepdistnew} holds with $\varphi_{p}\paren{r} = \min(1,cr^{-\alpha})$, and that it is exponentially weakly dependent if Equation~\eqref{eq:weakdepdistnew} holds with $\varphi_{p}\paren{r} = \min(1,c\exp\paren{-\gamma r})$, where $c\geq 1$, $\alpha>0$ (resp. $\gamma>0$) are some constants. In general, we are interested in the situation where
$\varphi_{p}(k) \rightarrow 0$ as $k\rightarrow \infty$, but where this decrease could be slow
(a form of ``long-range dependency'').

Under $(\WD(p))$, for $p\in [2,\infty)$ we obtain a non-asymptotic upper bound for $S_\cR$ in $L^p$-norm  (Burkholder-type inequality), and for $p=\infty$ we obtain that $S_\cR$ is sub-Gaussian with a bound on its sub-Gaussian norm (Azuma-Hoeffding type inequality).

We use the terminology of \textit{weakly dependent} random fields obeying Assumption~\ref{def:weak_spat_dep_new} by following the notions defined in the seminal works of \cite{Doukhan:94} and \cite{Bickel:99} and in the monograph \cite{Dedecker:07} which in a certain sense extend the concept of mixing to multi-parameter processes. 

The condition~\eqref{eq:weakdepdistnew} can be seen as a generalization in dimension $d$ of similar assumptions made
in dimension~1 in previous literature, and allowing the approximation of a sequence of random variables by a martingale difference,
also related to the notion of mixingale (see \cite{McLeish:76}). Various results based on assumptions
of a related nature are available in the literature, in particular for stationary sequences~\cite{peligrad:07,Peligrad:05} in dimension $1$,~\cite{Dedecker:15} also in dimension $1$ for Banach-valued processes, and~\cite{Dedecker:98,Dedecker:01}
for random processes in general dimension.

Specialized in dimension $d=1$, the condition~\eqref{eq:weakdepdistnew}
recovers similar conditions involving the $L_{p}$ norm of the projection of the past of random process, considered in~\cite{peligrad:07,Peligrad:05}.
In~\cite{Dedecker:98,Dedecker:01},
in general dimension $d$ the $\sigma$-algebra $\cM_{t,r}$ is intersected with the past $\sigma$-algebra of $t$ in the sense
of lexicographical ordering; this makes the condition in~\cite{Dedecker:98,Dedecker:01}
slightly weaker than $\WD(p)$; however
the lexicographical order seems to be chosen as an arbitrary total order over $\mbz^d$.
Furthermore, the 
condition~$\WD(p)$ 
implies bounds that improve over those of~\cite{Dedecker:01}
in the long-range dependency regime,
as will be discussed later.

In fact, the absence of a natural total order on $\mbz^d$ for $d>1$
is precisely what makes the generalization to higher dimension challenging and constitutes
the interest of our contribution; in a nutshell, we are able to generalize the dyadic martingale
decomposition used in~\cite{Peligrad:05,peligrad:07} to multi-scale decomposition in higher dimension by using a specific
total order on $\mbz^d$, and additional arguments.
For a comparison to the theory of strong multidimensional martingales, see Section~\ref{sec:disc}.


Since under our definition of polynomially/exponentially weakly-dependent processes it holds $\varphi_p(0)=1$, and since for any fixed $t$ the sequence $\norm{ \e{X_t | \cM_{t,k} } - \e{X_t}}_{p}$
is nonincreasing in $k\geq 0$ (by Jensen's inequality and the fact that the $\sigma$-algebras $\cM_{t,k}$ are
nested for fixed $t$), it is natural to think of $M_{p}$ as equal to $\sup_{t \in \mbz^d} \norm{X_{t} - \e{X_{t}}}_{p}$,
though we do not explicitly require it. Although we do not assume stationarity of the field, the fact
that $M_p$ is independent of $t$ indicates that our results are more geared towards fields that are in some sense
close to being stationary. However, motivated by applications in statistics, it is usual that instead of
a full random field over $\mbz^d$, one wants to consider only a process $(X_t)_{t\in \cD}$ defined over a subdomain $\cD$, typically a finite hyperrectangular window. Suppose that a version of Assumption $\WD(p)$ is satisfied for this process
but replacing $\mbz^d$ by $\cD$. Then if we extend $X_t$ to $\mbz^d$ by ``padding'' it with random variables that
are constant (say equal to zero) on $\mbz^d\setminus \cD$, we obtain a random field satisfying  Assumption $\WD(p)$ and can
apply our results (observe that the extended field is not stationary).

The paper is organized as follows. In Section~\ref{sec:main_res_01} we prove the main results.
In Section~\ref{sec:disc} we discuss the notion of multidimensional martingales and compare the results to the known bounds under different dependence measures. The multidimensional hierarchical martingale construction and supplementary technical lemmata are stated and proved in Section~\ref{sec:proof_main_bound} . 

\section{Main results}
\label{sec:main_res_01}


For any index set $\mathbb{T}$  and $\sigma-$fields $(\cF_{\alpha})_{\alpha \in \mathbb{T}}$ we write $\vee_{\alpha \in \mathbb{T}} \cF_{ \alpha }$ to denote the smallest $\sigma -$ field which contains all $\cF_{\alpha}$ for $\alpha \in \mathbb{T}$ and $\wedge_{\alpha \in \mathbb{T}}\cF_{ \alpha }$ for the $\sigma-$field which is the intersection of all $\sigma-$fields $\cF_{ \alpha }$. 

For a finite subset $\cD \subset \mbz^d$, we denote $S_{\cD} \bydef \sum_{t \in \cD}( X_{t} - \ee{}{X_{t}} )$ for the centered partial sum of $\paren{X_{t}}_{t \in \cD}$.
We denote $a\land b$, $a\lor b$ the minimum and maximum of two real numbers $a,b$.
We use the standard notation $\lceil a \rceil$, $[a]$,$\lfloor a \rfloor$ for the ceil, integer and floor part of a real number $a$. For any $k \in \mbn_{>0}$ we denote  $\irg{k}:= \{1,\ldots,k \}$.
We denote the cardinality of  a finite set $A$ as $\abs{A}$.
For a centered real random variable $Z$ we denote   $\norm{Z}_{\SG} \bydef \inf \set[1]{ c\geq 0: \e{\exp\paren{\lambda Z}} \leq \exp\paren[1]{\frac{\lambda^{2}c^{2}}{2}}, \text{ for all } \lambda \in \mbr } $ its subgaussian norm\footnote{It is well-known that
	the subgaussian norm is equivalent to the Orlicz $\psi_2$-norm for centered variables. We use the subgaussian norm here
	since it is the one appearing naturally in the Hoeffding-Azuma inequality.}. 




The following theorem constitutes our main contribution.

\begin{theorem}
	\label{thm:mres_rect2_new}
	Let $\cR=\prod_{i=1}^d \irg{N_i}$ be a $d$-dimensional rectangle of sidelengths $N_i \geq 1$, $i\in \irg{d} $, and $m(\cR):= \max_{i=1,\ldots,d} \lfloor \log_2 N_i \rfloor$.
	Let $\bmd=(\delta_k)_{k \geq 1}$ be a fixed {nondecreasing} sequence
	of nonnegative integers such that
	\begin{equation} 
		\label{eq:assptdel_01}
		\sum_{k=1}^{m(\cR)} \delta_k 2^{-k} \leq \frac{1}{4d^2}.
	\end{equation}

	Let $p\in[2,\infty]$ be fixed, and let $\paren{X_{t}}_{t \in \mbz^d}$
	be a random field such that Assumption~$\WD(p)$ is satisfied.
	\begin{itemize}
		\item if $p\in[2,\infty)$, then
		\[
		\norm{S_{\cR}}_{p} \leq 
  4\sqrt{p} \Psi_p(\bmd,\cR),
		\]
		\item if $p=\infty$, then
		\[
		\norm{S_{\cR}}_{\SG} \leq 10 \Psi_\infty(\bmd,\cR),
		\]
	\end{itemize}
	where (putting $\delta_0=0$)
	\begin{equation}
		\label{eq:psi_func}
		\Psi_{p}\paren{\bmd,\cR}  := 2M_p \sqrt{\abs{\cR}}\paren{2+ \sum_{k=1}^{m(\cR)+1}\varphi_{p}\paren[1]{{\delta}_{k-1}+1}\sqrt{\abs{\cC_{k} \cap \cR }}},
	\end{equation}
	with $\cC_{k} := \irg{2^k}^d$. 
\end{theorem}

Before proving the theorem, we give in the next corollary the rates obtained via the above bound for the partial sums of a weakly-dependent random field over cubes in $\mbz^d$  under the assumption that the random field $\paren{X_{t}}_{t \in \mbz^d}$ is either polynomially or exponentially weakly-dependent. 
In a nutshell, these are obtained
by choosing the sequence $\bmd=\paren{\delta_{k}}_{k\geq 1}$ such that the constraint~\eqref{eq:assptdel_01} is fulfilled
and the value of the function $\Psi_{p}\paren{\bmd,\cR}$ is close to its minimum.

\begin{corollary}
	\label{cor:specific_rates}
	Let $\cD_n=\irg{n}^{d}$ be a $d-$dimensional cube with side-length $n\geq 2$,
	and put $N=\abs{\cD_n}=n^d$.
	Consider a random field $\paren{X_{t}}_{ t \in \mbz^{d}}$ which satisfies the weak dependency Assumption~$\WD(p)$ for a given $p \in [2,+\infty]$ and the function $\varphi_{p}\paren{\cdot}$. Then 
	\begin{itemize}
		\item If $p \in [2,+\infty)$ and $\varphi_{p}\paren{r} \leq c r^{-\alpha}$, $\alpha>0$, $c\geq 1$,
		\begin{align}
			\label{eq:prate_cube}
			\norm[2]{N^{-1} {S_{\cD_n}}}_{p} & \leq c C_{d,\alpha} \sqrt{p} M_{p}N^{- \paren{ \frac{1}{2} \wedge \frac{\alpha}{d}}} , &\text{ if } \alpha \neq d/2 ;\\
			\norm{N^{-1} {S_{\cD_n}}}_{p} & \leq c \wt{C}_d \sqrt{p} M_{p}{N^{-\frac{1}{2}}} \paren{\log_{2}N}^{\frac{d}{2}+1}, &\text{ if } \alpha = d/2,
		\end{align}
		where $ C_{d,\alpha}:=24  \cdot 2^{d+\alpha} d^{2\alpha} \abs[2]{1-2^{\frac{\frac{d}{2} - \alpha}{1+\alpha}}}^{-(1+\alpha)}$, and
		$\wt{C}_d:= 15\cdot (8d)^{d/2}.$
		\item If $p =\infty$ and $\varphi_{\infty}(r) \leq ct^{-\alpha}$, then 
		\begin{align}
			\label{eq:irate_cube}
			\norm{N^{-1} S_{\cD_n}}_{\SG}& \leq 2.5 c C_{d,\alpha} M_{\infty}N^{- \paren{ \frac{1}{2} \wedge \frac{\alpha}{d}}}, & \text{ if } \alpha \neq d/2; \\
			\norm{N^{-1}S_{\cD_n}}_{\SG} & \leq 2.5 c \wt{C}_d M_{\infty} {N}^{-\frac{1}{2}} \paren{\log_{2}N}^{\frac{d}{2}+1}, & \text{ if } \alpha = d/2.
		\end{align}
	\end{itemize}
	Furthermore, if $\varphi_p\paren{t} \leq c \exp\paren{-\gamma t^{\eta}}$, then it holds $\norm{N^{-1}S_{\cD}}_{p} \leq c M_{p}C_{\gamma,\eta,d}{N^{-\frac{1}{2}}}$  \\ 
 ( $\norm{N^{-1}S_{\cD}}_{\SG} \leq c M_{\infty}C_{\gamma,\eta,d}{N^{-\frac{1}{2}}}$ when $p=\infty$ correspondingly ) for some factor $C_{\gamma,\eta,d}>0$ only depending on $\gamma,\eta,d$.
	
\end{corollary}

\begin{proof}
	{Let $\varphi_{p}\paren{t}\leq ct^{-\alpha}$.} Applying Theorem~\ref{thm:mres_rect2_new} to the cube $\cD = \irg{n}^{d}$,
	putting $m:=\lfloor \log_{2}n\rfloor$ we obtain that $\norm{S_{\cD}}_p$ resp. $\norm{S_{\cD}}_{\SG}$ is upper bounded proportionally to
	\[
	\Psi_{p}\paren{\bmd,\cD} \leq 2c  M_{p}\sqrt{N}\paren{2 + \sum_{k=1}^{m+1} \paren{\delta_{k-1}+1}^{-\alpha}\sqrt{2^{kd}} }, \]
	provided $\bmd$ satisfies the constraint~\eqref{eq:assptdel_01}. 
	
	Using Lagrange multiplier method 
	to minimize $\Psi_{p}(\bmd,\cD)$ under the constraint over $\bmd$ yields
	\[
	\wt{\delta}_{k} = \frac{1}{4d^2} Z_{\alpha,m}^{-1} (2\rho)^k 
	\qquad \text{ with }
	Z_{\alpha,m} := \sum_{k=1}^m \rho^k, \qquad \rho :=  2^{\frac{\frac{d}{2} - \alpha}{1+\alpha}}.
	\]
 {It is easy to check that $(2\rho) \geq 1$ so that
 $(\wt{\delta}_{k})$ is a nondecreasing sequence.}
	Take $\delta_{k} \bydef \big\lfloor \wt{\delta}_{k} \big\rfloor$, we have $\delta_k \leq \wt{\delta}_k$ so
	that constraint \eqref{eq:assptdel_01} is satisfied, {also
 $(\delta_k)$ is nondecreasing as required. Moreover,} $\delta_k +1 \geq \wt{\delta}_k$
	so that (recall $\delta_0=0$):
	\begin{align*}
		\Psi_{p}\paren{\bmd,\cD} \leq 2c M_{p}\sqrt{N}\paren{2 + 2^{\frac{d}{2}}\paren{1 +
				\sum_{k=1}^m 2^{\frac{kd}{2}} \wt{\delta}_k^{-\alpha}}}
		\leq 2c M_{p}\sqrt{N}\paren{2 + 2^{\frac{d}{2}}\paren{1 + 
				\paren{4d^{2}}^{\alpha} Z_{\alpha,m}^{1+\alpha}  }}.
	\end{align*}
	If $\alpha>\frac{d}{2}$ then $Z_{\alpha,m} \leq Z_{\alpha,\infty} = \rho(1-\rho)^{-1}$ 
	and we get
	\[\Psi_{p}\paren{\bmd,\cD} \leq \frac{1}{4} c M_{p}C_{\alpha,d}\sqrt{N} , \qquad
	C_{\alpha,d} := 24\cdot 2^{d+\alpha} d^{2\alpha}
	~
	\abs{1-\rho}^{-(1+\alpha)}.
	\]
	If $\alpha < \frac{d}{2}$ then $Z_{\alpha,m} \leq 
	\rho^{m+1}(\rho-1)^{-1} \leq
	N^{\frac{\frac{1}{2}-\frac{\alpha}{d}}{1+\alpha}}\rho \paren{\rho-1}^{-1}$  
	and
	\[
	\Psi_{p}\paren{\bmd,\cD} \leq \frac{1}{4} c M_{p}C_{\alpha,d} N^{1-\frac{\alpha}{d}} .
	\]
	Lastly, in the case $\alpha=\frac{d}{2}$ it holds $Z_{\alpha,m} =m \leq  d^{-1} \log_2 N $ and we get 
	\[
	\Psi_{p}\paren{\bmd,\cD} \leq 
	0.25 \cdot c M_{p} \wt{C}_{d} \sqrt{N} (\log N)^{\frac{d}{2}+1}, \qquad \wt{C}_{d} :=
	15 \cdot (8d)^{\frac{d}{2}}.
		\]
  {For  the case $\varphi_{p}\paren{t} \leq c \exp(- \gamma t^{\eta})$, it suffices to apply Theorem \ref{thm:mres_rect2_new} directly with $\delta_{k} = \big\lfloor \tilde{\delta}_{k}\big\rfloor$, $\tilde{\delta}_{k} = {2^{k/2}\paren[1]{\sqrt{2}-1}}/{4d^2}$. Note that such choice of $\delta_{k}$ ensures \eqref{eq:assptdel_01}. Thus we obtain that: 
  \[
	\Psi_{p}\paren[1]{\bmd,\cD} \leq 2c  M_{p}\sqrt{N}\paren{2 + \sum_{k=1}^{m+1} \exp\paren[2]{-C_{\gamma,\eta,d} \cdot 2^{\frac{k\eta}{2}}}{2^{\frac{kd}{2}}} }\leq cM_{p} C^{1}_{\gamma,\eta,d}\sqrt{N}, \]
 where $C_{\gamma,\eta,d}>0$, $C^{1}_{\gamma,\eta,d} >0$ some finite factors and the last inequality is due to \hfill \break
 $\sum_{k=1}^{\infty} \exp\paren[1]{-C_{\gamma,\eta,d} \cdot 2^{\frac{k\eta}{2}}}{2^{\frac{kd}{2}}} < \infty $.}

Finally the calculations for the case $p=\infty$ are identical and the claim follows.
	\end{proof}
\begin{remark}[On the optimality of the bound from Corollary \ref{cor:specific_rates}]
\label{rmk:opt_bound_spec}
{
In the case $\varphi_{p}(r) = \min(1,cr^{-\alpha})$ with $\alpha>\frac{d}{2}$ or $\varphi_{p}(r) = \min(1,c\exp(-\gamma r))$ ($\gamma >0$) the bounds \eqref{eq:prate_cube} and \eqref{eq:irate_cube} yield asymptotic behaviour (as $N$ grows) of order $\norm{N^{-1}S_{\cD}}_{p}= \cO\paren[1]{N^{-\frac{1}{2}}}$ and $\norm{N^{-1}S_{\cD}}_{SG} = \cO\paren[1]{N^{-1/2}}$ respectively, which is known to be optimal if the field
is i.i.d.
    Let us consider the case of polynomially weakly dependent processes with rate $\varphi_{p}(r) = \min(1,cr^{-\alpha})$ and $\alpha<\frac{d}{2}$ (we will not comment
    here on the optimality of the additional $\log(N)$ factor
    if $\alpha = \frac{d}{2}$).
 Consider the field $X_{t}\paren{\omega} := m_0 \xi \mbi_{t \in \irg{n}^d}$, where $m_{0}=n^{-\alpha}$ and $\xi$ is a Rademacher variable. In other words $X_t$ is constant
 equal to $\pm m_0$ with a randomly flipped sign over the cube of side-length $n$ and 0 elsewhere. One can readily check that $(X_{t})_{t \in \mbz^d}$ satisfies  $\WD(p)$ with $M_{p} = 1$ and furthermore $\abs{ N^{-1}S_{\cD}} = m_0 =N^{-\frac{\alpha}{d}}$. This example shows that inequalities~\eqref{eq:prate_cube},\eqref{eq:irate_cube}
 are tight under the given assumptions, up to 
 multiplicative factors only depending on $p$ and $d$.
In other words the bounds of Corollary~\ref{cor:specific_rates}
are saturated by either of the extremal cases of an i.i.d. field 
or a constant field on the cube. Concerning the latter,
note that the exhibited distribution saturating the bound depends on the side-length $n$.
It would be more satisfactory to exhibit an example
of a {\em fixed} field distribution satisfying
the polynomial weak mixing condition with $\alpha < \frac{d}{2}$
("long range" or "heavy tail" dependence) and such
that the rates~\eqref{eq:prate_cube},\eqref{eq:irate_cube}
are tight asymptotically as $N$ grows.
}
\end{remark}
 
	\section{Discussion and comparison to known results}
	\label{sec:disc}
	\newcommand{\bn}{{\mbf{n}}}
	In this section we compare our approach and bounds 
    to some existing results for
	dependent real random fields. In many of such results (see ex. \citep{Doukhan:84,Rio:00,Dedecker:01}) an exponential inequality for partial sums of bounded random fields is derived from $L_{p}\paren{\mbp}$ bound by optimizing over the value $p\geq 2$ . {The  analysis of the properties of dependent random fields are 
    commonly based  on approximations by 
    multidimensional martingales.
    For this reason we first review briefly the different notions of multidimensional martingales, associated Burkholder-type inequalities for this reference case, and the relation to approximations of dependent random fields. In the following section we also discuss the case of Bernoulli random fields, that is, stationary transforms of an i.i.d. random field. }
	
	Notice that the weak dependency condition~\eqref{eq:weakdepdistnew} when $d=1$ is stronger than
 the so-called "mixingale type" condition (see for example \cite{McLeish:76}), for which the process is not necessarily adapted to the conditioning filtration. The former is mentioned in the works \cite{Dedecker:07},\cite{Dehling:82} to
	characterize fading correlation between the past and the future of a discrete stochastic process. 
	\subsection{Multiparameter martingale difference fields} 
\label{sec:multip_mart_diff}
{Although our interest concerns chiefly dependent random fields, we first review briefly
several notions of multiparameter martingale difference fields as reference situations.}
For $k,\ell \in \mbz^{d}$ we denote $k \leq_{cw} \ell$ if the inequality holds coordinate-wise, i.e. $k_{i}\leq \ell_{i}$ for all $i \in \irg{d}$, $k <_{cw} \ell$ if the inequality is strict for at least one coordinate, and $k \ll_{cw} \ell$ if the inequality is strict for {\em all} coordinates. Simplifying the general point of view to focus on
canonical filtrations in our setting,
define the $d$ {\em marginal} filtrations $(\cF^{(k)}_{i} := \fS \set{X_{(t_1,\ldots,t_d)}: t_k \leq i})_{i \in \mbz}$, $k\in \irg{d}$; and for $t=(t_1,\ldots,t_d)\in \mbz^d$, let $\cF^\wedge_t := \bigwedge_{k \in \irg{d}}\cF^{(k)}_{t_k} $, and $\cF^\vee_t := \bigvee_{k \in \irg{d}}\cF^{(k)}_{t_k} $.
Note that  $\cF^\vee_t  = \cM_{t,0}$.
The field $(X_t)_{t\in\mbz^d}$ is then called a
\begin{enumerate}[label=(\roman*)]
\item {\em weak} martingale difference field, if $\e{X_t | \cF^{\wedge}_{t'}}=0$ for all $t' <_{cw} t$;
\item {\em ortho-martingale} difference field, if $\e[1]{X_t | \cF^{(k)}_{t'_k}}=0$ for all $t=(t_1,\ldots,t_d) \in \mbz^d$, $k\in \irg{d}$
and $t'_k<t_k$;
\item {\em strong} martingale difference field, if $\e{X_t | \cF^{\vee}_{t'}}=0$ for all $t' \ll_{cw} t$.
\end{enumerate}

The above terminology is in line with the general definitions introduced in the seminal paper \cite{Cairoli:75}, see also \cite{Basu:79}.
Since $ \cF^{\wedge}_{t'} \subseteq \cF^{(k)}_{t'_k} \subseteq \cF^{\vee}_{t'}$,
we have the straightforward implication $(iii) \Rightarrow (ii) \Rightarrow (i)$.
It is also well-known that $(i)$ implies $(ii)$ under the additional assumption of {\em commuting}
marginal $\sigma$-fields
also known as condition (F4) from \cite{Cairoli:75}, see
\cite{Khoshnevisan:02} for a modern account and equivalent formulations
of the assumption. {However, we stress that $(iii)$ does not imply
commuting marginal $\sigma$-fields in general.}

{The projective-type dependency assumption $\WD(p)$, in the edge case $\varphi_{p}(0)=1$ and $\varphi_{p}(t)= 0$ for $t>0$, reduces to the strong martingale difference field assumption $(iii)$. On the other hand, under assumption $(iii)$ a Burkholder-type inequality for random fields can be derived as a straightforward consequence of Burkholder's inequality in dimension $1$.} This holds since $(iii)$ implies a (one-dimensional) martingale difference condition under the filtration $\cF_{t}^{\lelex} := \fS\{ X_{t'}: t' \lelex t \}$ and the one-dimensional Burkholder's inequality can be applied to the elements of the partial sum ordered in lexicographic order, i.e, we obtain 
\begin{equation} \label{eq:compar} 
\norm[1]{S_{\irg{n}^d}}_{p} \leq C_{p}n^{\frac{d}{2}}\max_{t \in \irg{n}^d}\norm{X_t}_{p}
\end{equation} 
by this direct argument.  {Therefore, 
in the particular case $(iii)$ our results do not bring something new (note also that in that case the optimal
multiplicative constant 
in the bound does not depend on the dimension $d$, while ours does).
Still, from this point of view the assumption $\WD(p)$ in the general case can be interpreted as a form of limited departure of the random field from a strong martingale difference field.}

 {Concerning the other types of multi-parameter martingales, }
 in the bivariate case a Burkholder-type inequality was established
 by~\cite{Gundy:76} (Lemma 1 there) under assumptions $(i)$+(F4), and under assumption $(ii)$ by~\cite{Metraux:78}.
 For general dimension $d$, \cite{Faz:05} established
 Burkholder-type inequalities, again under the assumptions $(i)$+(F4). Both  \cite{Gundy:76} and \cite{Faz:05} use an iterative approach over 
 the dimension. These results are not recovered by our approach since the 
 assumptions are different.

 {Recent developments have used notions of approximation of a random field by ortho-martingales to obtain various results of a related nature
\citep{Volny:15,Elmachkouri:16,Giraudo:18,Giraudo:19}.
 However, we note that these works use the terminology ``ortho-martingale'' to mean in fact $(i)$+(F4). We insist that we do not 
 use the assumption (F4) of commuting marginal filtrations in the present work,
 so that we see our results as going in a different direction in nature.}

\subsection{ Bernoulli random fields}
\label{sec:bern_rand_fields}
A Burkholder-type inequality for so-called Bernoulli random fields is given in \cite{Wu:13}. In that work a class of random fields of the form  $X_{t} = g\paren{\epsilon_{t-s},s \in \mbz^d}$, $t \in \mbz^d$ where $\paren{\epsilon_{i}}_{i \in \mbz^d}$ are i.i.d. random variables and $g$ is some bounded measurable function is considered. Typical examples which belongs to this class are linear random fields and Volterra random fields (see examples in \cite{Sang:17} and Section~2 in \cite{Giraudo:19}).  Functional central limit theorems \cite{Bierme:14,Klicnarova:16}, a large deviation principle \cite{Sang:17} and  a variant of law of the iterated logarithm for partial sums of this class of random fields \cite{Giraudo:22} have been established.  For these types of random fields a following dependence measure (originally introduced in \cite{Wu:05}) is introduced. Namely, for a stochastic random field $\paren{X_{t}}_{t \in \mbz^d}$, such that $X_{t} \in L_{p}\paren{\mbp}$ define the coefficient $\delta_{t,p}=\norm{X_t -X_{t}^{\star}}_{p}$, where the coupled process $X^*_t$ is given by $X^{\star}_{t}=g\paren{\epsilon^{\star}_{t-s}, s \in \mbz^{d}}$ and $\epsilon^{\star}_{j} = \epsilon_{j}\mbi_{j \neq \bm{0}}+\epsilon'_{\bm{0}}\mbi_{j = \bm{0}}$, with 
$\epsilon'_{\bm{0}}$ being 
{an independent copy of $\epsilon_{\bm{0}}$.}
In this framework, Proposition~1 of \cite{Wu:13} implies following the Burkholder-type inequality for the partial sums of $X_{t}$:
\begin{align*}
\norm[1]{S_{\irg{n}^d}}_{p} \leq \sqrt{2p}{n^{d/2}}\Delta_{p},
\end{align*}  
where $\Delta_{p} = \sum_{i \in \mbz^d}\delta_{i,p} < \infty$.

{In the one-dimensional case,
a particular situation is when the field $(X_t)_{t \in \mbz}$ is {\em causal} with respect to the underlying i.i.d. field $(\eps_t)_{t \in \mbz}$, that is,
$X_t=g(\eps_{t-s}, s\geq 0)$. In this
situation, as pointed out in \cite{peligrad:07}, if we denote
$\cF_k:=\sigma((\eps_t)_{t\leq k})$, it holds by Jensen's  and stationarity
\begin{align*}
\norm{\e{X_t|\cM_{t,k} }- \e{X_t}}_p
& = \norm{\e{X_k|\cM_{k,k}} - \e{X_k}}_p\\
& \leq \norm{\e{X_k|\cF_{0}} - \e{X_k}}_p\\
&\leq \norm{g(\eps_k,\eps_{k-1},\ldots)
-g(\eps_k,\eps_{k-1},\ldots,\eps_0,\eps'_{-1},\eps'_{-2},\ldots)}_p,
\end{align*}
so that assumption \WD\ can be related to quantities that are tractable in
many models (see \cite{peligrad:07} for more discussion).}  {Notice that an analogous notion of a "causal field" in higher dimension would be
of the form $X_t=g(\eps_u, u \not \gg_{cw} t)$
(with the notation $\gg_{cw}$ as introduced in Section~\ref{sec:multip_mart_diff}), a notion that has been considered in past literature (e.g. \cite{Giraudo:21b}, Corollary 3.5). Observe that in the Bernoulli random field model, the marginal filtrations
for the i.i.d. field $(\eps_t)_{t \in \mbz^d}$ are commuting, so that tools based on ortho-martingale approximation mentioned in the previous
section can be applied as well (see ex. \cite{Giraudo:21,Giraudo:21b}).
In such a setting we can, by the same token as in
dimension 1, relate assumption \WD\ to tractable quantities.
To deal with non-causal Bernoulli fields however will require to relax assumption \WD\ to the non-adapted case (wherein we condition with respect
 to the field $(\eps_t)_{t \in \mbz^{d}}$ rather than $(X_t)_{t \in \mbz^d}$);
which is out of the scope of the present work.}

\subsection{{$L_p(\mbp)-$ projective criterion}}

In the work \cite{Dedecker:01}, a Burkholder-type inequality for random fields is obtained under an $L_{p}-$projective dependence criterion; it takes the form
\begin{align}
	\label{eq:proj_depen}
	\norm[1]{ S_{ \irg{n}^d}}_{p} \leq \sqrt{2p \sum_{t \in \irg{n}^d}b_{t,p/2}\paren{X}},
\end{align}
where $b_{t,\alpha}(X) = \norm{X_{t}^{2}}_{\alpha} + \sum_{k \in V_{t}^{1}}\norm{X_{k} \ee{\abs{k-t}}{X_{t}}}_{\alpha}$ (we assume the
field is centered for simplicity),
$p\geq 2$, $\abs{x-y}= \max_{1\leq i \leq d}\abs{x_{i}-y_{i}}$, $V_{t}^{1}$ denotes the set of all elements which precede $t$ in lexicographic order on $\mbz^d$ and $\ee{\ell}{X_{t}}$ is the conditional expectation with respect to the $\sigma-$algebra $\cF^{}_{V^{\ell}_{t}}$, where $V_{t}^{\ell}$ is the subset of elements from $V^{1}_{t}$ which are at distance at least $\ell$ from $t$. 
In the case where $\paren{X}_{t \in \cD}$ is a strong martingale difference random field then $b_{t,p/2} = \norm{X_{t}^{2}}_{p/2} = \norm{X_{t}}_{p}^{2}$ and \eqref{eq:proj_depen} implies \eqref{eq:compar} with a constant $C_{p,d} = \cO(\sqrt{p})$. 

{Let us analyze briefly the estimates obtained via bound~\eqref{eq:proj_depen} under the projective weak
dependency assumption $\WD'(p)$
(that is, $\WD(p)$ with $\cM'_{t,r}$ in place of
$\cM_{t,r}$, see \eqref{eq:mprime}).}
By construction of the $\sigma-$ algebras, it holds $\cF^{}_{V_{t}^{r}} \subseteq {\cM'_{t,r}}$, 
{thus under assumption $\WD'(p)$, by
the generalized H\"older's inequality
it holds $\norm{X_k \ee{|k-t|}{X_t}}_{p/2} \leq \norm[1]{X_k}_p \norm[1]{\e[1]{X_t|\cM'_{t,|k-t|}}}_p
\leq M_p^2\varphi_p(|k-t|)$.}
Hence,
for the right-hand-side of Inequality \eqref{eq:proj_depen} when $p \geq 2$ we get
\begin{align*}
	\sqrt{2p \sum_{t \in \irg{n}^d}b_{t,p/2}\paren{X}} 
	& \leq {M_p\sqrt{2p \paren[2]{\sum_{t \in \irg{n}^d } \paren[2]{1 + \sum_{r=1}^{\infty}r^{d-1}\varphi_{p}(r) }  } },}
\end{align*}
Now, if $\varphi_{p}(r)\lesssim r^{-\alpha}$, 
{plugging this into \eqref{eq:proj_depen} yields
a bound of order $O\paren[1]{n^{-\paren[1]{\frac{d}{2}\wedge\frac{\alpha}{2}}}}$ 
for $n^{-d}\norm[1]{ S_{ \irg{n}^d}}_{p}$, while Corollary \ref{cor:specific_rates} entails 
a bound in~$O\paren[1]{n^{-\paren[1]{\frac{d}{2} \wedge \alpha}}}$.
Thus, under assumption \WD' the estimate obtained from Corollary~\ref{cor:specific_rates}
improves over that obtained via the bound~\eqref{eq:proj_depen}
(in particular extending the range of $\alpha$ for which a convergence of
order $O(n^{-d/2})$ is granted).
On the other hand, it should be noted that under the stronger 
assumption of control of the $\alpha$-mixing coefficient
between the $\sigma$-algebras $\fS(X_t)$ and $\cM'_{t,k}$, using Rio's covariance inequalities \cite{Rio:93} one can infer both assumption $\WD'(p)$ for a certain
function $\varphi_p(k)$ but also
$\norm{X_k\ee{|k-t|}{X_t}}_p/2 \lesssim \varphi_p(|k-t|)^2$. In this scenario, if $\varphi_p(r) \lesssim r^{-\alpha}$, 
bound~\eqref{eq:proj_depen} and Corollary~\ref{cor:specific_rates} yield a bound of the same order \cite{Dedecker:20}.}

\section{Proof of main result}
\label{sec:proof_main_bound}
To prove Theorem~\ref{thm:mres_rect2_new} it turns out to be more convenient to consider rectangles starting at the origin. We therefore introduce the notation
$\irg{k}_{0} \bydef \{0,\ldots,k-1 \}$ if $k \geq 1$ and $\irg{k}_{0} = \emptyset$ if $k=0$, and
below we consider rectangles of the form $\cR=\prod_{i=1}^d \irg{N_i}_0$.

For sets $A,B \subseteq \mbz^d$, and $c \in \mbz$, $v \in \mbz^{d}$ we will use the standard notation
\begin{align*}
	c A & := \set{ca: a \in A} \subseteq \mbz^d;\\
	A+B & := \set{a+b: a \in A, b \in B}  \subseteq \mbz^d.
\end{align*}
(If $A=\emptyset$ we denote $ cA := \emptyset$, 
$A+B :=\emptyset$.) 

{
Before we begin the proof proper, we give a few words
of informal overview of its structure. The main principle
of the proof is a decomposition of each term $X_t - \e{X_t}$
as a "telescopic" series with terms $\e{X_t|\cF^{\prec}_k(t)} - \e{X_t|\cF^{\prec}_{k-1}(t)}$, where $\cF^{\prec}_k(t)$
is a suitable filtration generated by blocks 
(which we will call {\em cells}) of the 
random field at a dyadic scale $2^k$. Each of these terms
can be then summed over the corresponding dyadic cells and
gives rise to a martingale difference "at scale $k$" (i.e.
indexed by dyadic integer vectors);
for each $k$ the corresponding dyadic martingale can be controlled by the Marcinkiewicz-Zygmund inequality.
This canvas follows the general line of arguments
used by~\cite{Peligrad:05,peligrad:07} in dimension $d=1$.}

{Going to higher dimension, we run into the following
additional challenges:
\begin{itemize}
\item  There is no natural total order on $\mbz^d$. In order
to use martingale arguments, we need to carefully define
a specific total order on $\mbz^d$ to get the martingale structure
while ensuring a form of compatibility of filtrations across scales; 
this construction will be explained in detail in Section~\ref{subsec:MMD}.
\item We need to control
$\norm{\e{X_t|\cF^{\prec}_k(t)} - \e{X_t|\cF^{\prec}_{k-1}(t)}}_p$
using the $\WD(p)$ assumption,
and this is directly linked
to $\varphi_p(d_k(t))$, where
$d_k(t)$ is the (supremum) distance
of $t$ to the border of the dyadic cell it belongs to. Directly summing 
these quantities over cells gives rise to good estimates in dimension 1
(see~\cite{peligrad:07}), but it turns out that the estimates obtained this
way for rectangles in dimension $d\geq 2$ are suboptimal.
The reason is that there are ``too many'' elements
in the cell that are close to its boundary,
making a too large contribution in the
sum over the cell.
This is why, to alleviate this issue, we will at first 
exclude the elements which are close to boundaries of a cell at any scale, which we will call the (multiscale) "frame", so that the remaining elements are then sufficiently ``separated'' from the cell boundaries. The sum of elements in the "frame" will then be dealt with by an inductive argument.
This framed decomposition and the inductive argument will be explained
first, in the coming section.
\end{itemize}
}

\subsection{Framed decomposition and main inductive step}

\label{subsec:FDIS}

We prove the theorem by induction on the size of the rectangle $\cR$.

We consider the case $p<\infty$ only; the arguments for the case $p=\infty$ are the same for the remainder of the proof.
For $\cR$ reduced to a single element $\bm{0}$, the claim obviously holds. Now, assume the claim is established for any rectangle $\cR' \subsetneq \cR$. 
We use the following construction: 
for any integer $\delta< 2^k$  let
\begin{equation}
	\label{eq:defLambda_01}
	\Lambda_{k,\delta} := 2^k\mbn_{>0} + \irg{\delta}_0. \qquad (\text{If } \delta=0, \text{ then } \Lambda_{k,\delta} := \emptyset.)
\end{equation}

\begin{figure}[h!]
	\vspace*{0.2cm}
	\centerline{\resizebox{0.6\linewidth}{!}{{\adjustbox{Clip*= 0cm 0cm 13cm 10cm}{\resizebox{1.5\linewidth}{!}{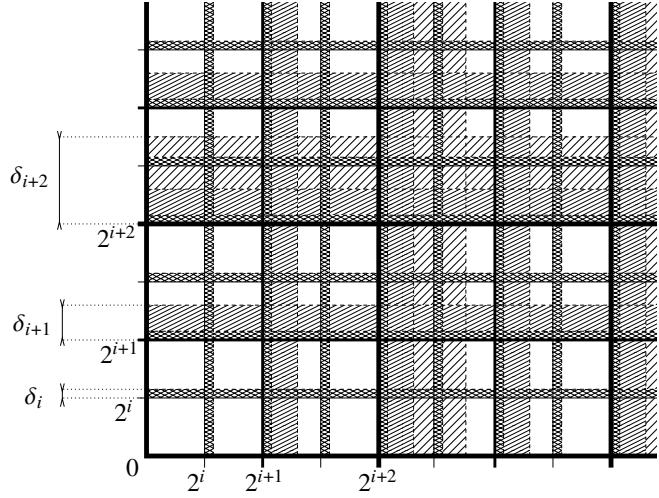}}}}}
	\caption{ \label{fig:figframe} Illustration of ``framed set'' and
		the ``frame'' for the case $d=2$ on three consecutive scales. \textit{Hatched regions} belong to the frame $\fF_\bmd^c$. The \textit{white region} is the framed set $\fF_\bmd$.}
\end{figure}

Let $\bmd=(\delta_k)_{k \geq 1}$ be a fixed sequence
of integers with $\delta_k \leq 2^k$, $k\geq 1$. Define
\begin{align} 
	\label{eq:diff_sets}
	\Lambda_{\bmd} := \bigcup_{k \geq 1} \Lambda_{k,\delta_k} \qquad
	\fF_{\bmd} :=  \paren{ \mbn\setminus\Lambda_{\bmd}}^d ; \qquad \fF^c_{\bmd} := \mbz^d \setminus \fF_{\bmd}.
\end{align}
We call $\fF_{\bmd}$ the ``framed set'' and $\fF^c_{\bmd}$ the ``frame''. These sets are illustrated
on Figure ~\ref{fig:figframe}.  
We use the decomposition $ \cR = (\cR \cap \fF_{\bmd}) \uplus (\cR \cap \fF^c_{\bmd})$. By the triangle inequality 
\begin{align}
	\label{eq:framedcmp_new}
	\norm{S_{\cR}}_{p} \leq \norm{S_{\cR \cap \fF_{\bmd} }}_{p} + \norm{S_{\cR \cap \fF^{c}_{\bmd} }}_{p}.
\end{align} 

The proof proceeds as follows. The first and main term will be controlled by a multiscale martingale
decomposition and lead to the crucial estimate $\norm{S_{\cR \cap \fF_{\bmd} }}_{p} \leq \frac{C_{p}}{2}\Psi_{p}\paren{\bmd,\cR}$ where
$C_p = 4 \sqrt{p}$. The proof of this estimate is postponed to
Proposition~\ref{prop:framedestimate_01} below.

The second, remainder term will be dealt with via the induction hypothesis. We first decompose $\fF_\bmd^c$ as a disjoint union of
product sets, writing $\Lambda_\bmd^c := \mbn \setminus \Lambda_\bmd$:
\[
\fF_\bmd^c = \mbn^d \setminus (\Lambda_\bmd^c)^d = \biguplus_{i=1}^d \Delta_i, \qquad \Delta_i := \paren{\prod_{j=1}^{i-1} \Lambda_{\bmd}^c \times \Lambda_{\bmd} \times  \prod_{j=i+1}^d \mbn },
\]
therefore, by the triangle inequality,
\begin{equation}
	\label{eq:sumcompr_01}
	\norm{ S_{\cR \cap \fF_\bmd^c}}_p = \norm[3]{\sum_{i=1}^d  S_{\cR\cap \Delta_i}}_p \leq \sum_{i=1}^d \norm{ S_{\cR\cap \Delta_i}}_p.
\end{equation}
We introduce the following notation. For a finite set $B \subset \mbn$, and an integer $j\in \irg{\abs{B}}_0$, denote
$(j: B)$ the $(j+1)$-th element of $B$ in increasing order. For a finite product set $\bA=\prod_{i=1}^d A_i \subseteq \mbn^d$,
define $\cK(\bA):=\prod_{i=1}^d \irg{\abs{A_i}}_0$, 
and for
a $d$-tuple ${\bm t}=(t_1,\ldots,t_d) \in \cK(\bA)$, denote $(t:\bA)=((t_1:A_1),\ldots,(t_d:A_d))$,
and the ``compressed'' version of the restriction of the process $(X_t)_{t \in \mbz^{d}}$ to $\bA$ as
\begin{align} 
	\label{def:compressed_set_def}
	\begin{cases}
		\wt{X}^{(\bA)}_t = X_{t:\bA}, &  t \in \cK(\bA),\\
		\wt{X}^{(\bA)}_t = 0, & t \not\in \cK(\bA);\\
	\end{cases}
	\quad \quad  \text{ then it holds } \qquad
	\wt{S}_{\cK(\bA)}^{(\bA)} := \sum_{t \in \cK(\bA)} \wt{X}^{(\bA)}_t = S_\bA. 
\end{align}

Since $\Delta_i$ is a product set, so is $\cR\cap \Delta_i$, and we can apply the above ''compression principle''.
Using assumption~\eqref{eq:assptdel_01}, the side-length of the rectangle $\cK(\cR\cap \Delta_i)$ along direction $i$ is bounded by
\[
\abs[1]{\Lambda_\bmd \cap \irg{N_i}_0} = \abs[4]{\bigcup_{k\geq 1} \Lambda_{k,\delta_k} \cap \irg{N_i}_0 }
\leq \sum_{k=1}^{\lfloor \log_2(N_i) \rfloor} \bigg\lfloor \frac{N_i}{2^k} \bigg\rfloor \delta_k
\leq N_i \sum_{k=1}^{m(\cR)} 2^{-k}\delta_k \leq \frac{N_i}{4d^2},
\]
while for $j \neq i$ the side-lengths are bounded by $N_{j}$. Therefore $\cK(\cR \cap \Delta_i) \subsetneq \cR$ and
\begin{equation}\label{eq:bndcard_01}
	\abs{\cR \cap \Delta_i} \leq \frac{\abs{\cR}}{4d^2}.
\end{equation}
By Lemma~\ref{lem:compression_property} below, the ``compressed'' process $\paren[1]{\wt{X}^{(\cR\cap\Delta_i)}_t}_{t \in \mbz^{d}}$ satisfies the same weak-dependency condition \eqref{eq:weakdepdistnew} as the original process.  Applying the induction hypothesis to the process $\paren[1]{\wt{X}^{(\cR\cap\Delta_i)}_t}_{t \in \mbz^{d}}$ over
the rectangle $\cK(\cR\cap\Delta_i)$, we obtain
\begin{equation} \label{eq:estcompr_01}
	\norm{S_{R\cap \Delta_i}}_p = \norm[1]{\wt{S}_{\cK(R\cap \Delta_i)}^{(R\cap \Delta_i)}}_p \leq C_{p} \Psi_p(\bmd,\cK(\cR \cap \Delta_i)).
\end{equation}
We estimate this upper bound using~\eqref{eq:bndcard_01} and straightforward cardinality bounds via:
\begin{align}
	\Psi_{p}(\bmd,\cK(\cR \cap \Delta_i))
	& =2M_p \sqrt{\abs{\cR \cap \Delta_i }}\paren[3]{2+\sum_{k=1}^{m(\cK(\cR \cap \Delta_i))+1}\varphi_{p}\paren[1]{{\delta}_{k-1}+1}\sqrt{\abs{\cC_{k,0} \cap \cK(\cR\cap\Delta_i) }}}
	\notag\\ 
	&\leq
	\frac{1}{d} M_p \sqrt{\abs{\cR}}\paren[3]{2+\sum_{k=1}^{m(\cR)+1}\varphi_{p}\paren[1]{{\delta}_{k-1}+1}\sqrt{\abs{\cC_{k,0} \cap \cR} }} \notag\\
	&=
	\frac{1}{2d}   \Psi_{p}(\bmd,\cR). \label{eq:estpsi_01}
\end{align}
Finally combining Equation~\eqref{eq:framedcmp_new}, Proposition~\ref{prop:framedestimate_01}, Equations~\eqref{eq:sumcompr_01},~\eqref{eq:estcompr_01} and~\eqref{eq:estpsi_01}, we obtain
\[
\norm{S_\cR}_p \leq \frac{1}{2} C_{p} \Psi(\bmd,\cR) + \sum_{i=1}^d \frac{1}{2d} C_{p} \Psi_p(\bmd,\cR) \leq C_{p} \Psi_p(\bmd,\cR),
\]
and the induction claim is proved. \qed

\begin{lemma}
	\label{lem:compression_property}
	Let $p \in [2, \infty]$ be fixed and assume  the process $\paren{X_{t}}_{t \in \mbn^{d}}$ satisfies
	assumption $\WD(p)$ with $\varphi_{p}\paren{\cdot}$. Then for any product set $\bA \subset \mbn^{d}$,
	the "compressed" version $\paren[1]{\wt{X}^{(\bA)}_{t}}_{t \in \mbn^{d}}$ as defined by \eqref{def:compressed_set_def} satisfies
	assumption $\WD(p)$ with the same function $\varphi_{p}\paren{\cdot}$.
\end{lemma}
\begin{proof}
	For every $t \notin \cK\paren{\bA}$, we have by construction that $\wt{X}^{\paren{\bA}}_{t} =0$,
	it is therefore sufficient to establish property~\eqref{eq:weakdepdistnew} only for elements $t \in \cK\paren{\bA}$.
	Recall $\cM_{t,k} = \fS \{ X_{u}: u \in \mbz^d, 
     \sup_{i}(t_i - u_i) \geq k
     \}$ and let
	\begin{align*}
		\wt{\cM}_{t,k}^{(\mbf{A})} & \bydef \fS \{ \wt{X}^{(\mbf A)}_{u}: u \in \mbz^d, 
        \sup_{i}\paren{t_{i} - u_i} \geq k
        \}\\
		& = \fS \{ X_{(u:\mbf{A})}: u \in \cK(\mbf{A}), 
        \sup_{i}\paren{t_{i} - u_i} \geq k
        \},
	\end{align*}
	wherein we are able to restrict for $u \in \cK(\mbf{A})$ using again $\wt{X}^{\paren{\bA}}_{u} =0$ for $u \not\in \cK(\mbf{A})$.
	
	For every $t,u \in \cK\paren{\bA}$, it holds 
	$\sup_{i}\paren{ (t_{i}:A_{i}) - (u_{i}:A_{i})} \geq \sup_{i}\paren{t_{i} - u_i}$, so for any $k>0$: 
	\begin{align*}
		\wt{\cM}_{t,k}^{(\mbf{A})} = \fS \{ X_{(u:\mbf{A})}: u \in \cK(\mbf{A}), 
        \sup_{i}\paren{t_{i} - u_{i}}\}
		& \subseteq \fS \{ X_{(u:\mbf{A})}: u \in \cK(\mbf{A}), 
        \sup_{i}\paren{t_{i}:A_{i} - u_i : A_{i}} \geq k \}\\
		& \subseteq \fS \{ X_{v}:  
        \sup_{i}\paren{t_{i}:A_{i} - v_{i}}\geq k \}\\
		& = \cM_{(t:\bA),k}.
	\end{align*}
By using Jensen's inequality, for any $t\in \cK(\bA)$, $k> 0$, using the assumption $\WD(p)$ for process $(X_t)$ we get: 
	\begin{align*}
		\norm{\e{\wt{X}^{\paren{\bA}}_{t}|\wt{\cM}^{(\bA)}_{t,k}} -\e{\wt{X}^{\paren{\bA}}_{t}} }_{p}
		& =     \norm{\e{X_{(t:\bA)}|\wt{\cM}^{(\bA)}_{t,k}} -\e{X_{(t:\bA)}}}_{p}\\
		& \leq  \norm{\e{X_{(t:\bA)}|\cM_{(t:\bA),k}} -\e{X_{(t:\bA)}}}_{p} \leq M_p \varphi_p(k),
	\end{align*}
	so process $(X^{(\bA)}_t)$ satisfies $\WD(p)$.
\end{proof}

\subsection{Multiscale martingale decomposition}
\label{subsec:MMD}

As announced previously, the following estimate for $\norm{S_{\cR \cap \fF_{\bmd}}}_p$ is crucial.
The main argument resides on a multi-scale martingale decomposition, the principle of which was
used in dimension $d=1$ in \cite{peligrad:07}. However, in dimension $d\geq 2$
the conditioning $\sigma$-algebra
$\cM_{t,k}$ in Assumption $\WD(p)$ cannot be expressed as a ``past'' for a total order on $\mbz^d$,
and extending the argument requires a more involved construction.

\begin{proposition}
	\label{prop:framedestimate_01}
	Let $\cR=\prod_{i=1}^d \irg{N_i}_0$ be a $d$-dimensional rectangle of side-lengths $N_i \geq 1$, $i=1,\ldots,d$, and $m(\cR):= \max_{i=1,\ldots,d} \lfloor \log_2 N_i \rfloor$.
	Let $\bmd=(\delta_k)_{k \geq 1}$ be a fixed {nondecreasing} sequence
	of integers with $\delta_k \leq 2^k$, $k\geq 1$ (put $\delta_0=0$), and let $\fF_{\bmd} $ be as defined in \eqref{eq:diff_sets}.

	Assume the process $\paren{X_{t}}_{t \in \mbn^d}$ satisfies the weak dependency assumption
	$\WD(p)$ for some with $p \in [2, \infty]$. 
	If $p<\infty$, it holds
	\begin{equation}
		\label{eq:contframedp_01}
		\norm{S_{\cR \cap \fF_{\bmd}}}_p \leq \frac{1}{2} C_{p} \Psi_p(\bmd,\cR),
	\end{equation}
	where $C_{p}:=4 \sqrt{p}$;
	and if $p=\infty$,
	then it holds
	\begin{equation}
		\label{eq:contrframedsubg}
		\norm{S_{\cR \cap \fF_{\bmd}}}_{\SG} \leq \frac{1}{2} C_{\infty} \Psi_p(\bmd,\cR),
	\end{equation}
	where $C_{\infty} := 10$ and we recall that $\Psi_{p}(\bmd,\cR)$ is given by \eqref{eq:psi_func}.
\end{proposition}

\begin{figure}[h!]
	\begin{center}
		\resizebox{.8\linewidth}{!}{\centerline{\adjustbox{Clip*= 0cm 0cm 13cm 10cm}{\resizebox{.8\linewidth}{!}{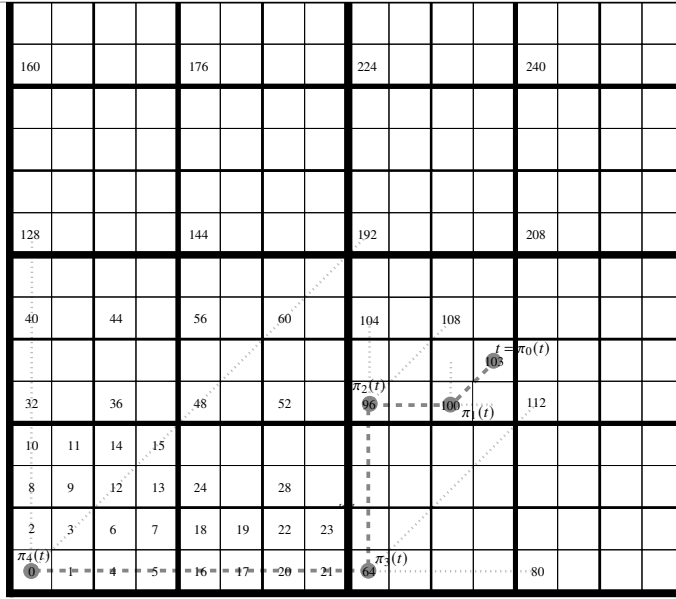}}}}
	\end{center}
	\caption{Illustration of the total order $\preceq$ for $d=2$. The numbers indicate the
		total ordering of the cells of $\mbn^2$ according to $\preceq$.
		The successive ``dyadic projections'' $\pi_k(t)$ for a particular point
		$t$ are given (dashed line). For each projection $\pi_k(t)$, dotted lines
		point to the possible positions of $\pi_{k-1}(t')$ such that $\pi_k(t)=\pi_k(t')$.}
	\label{pic:order}
\end{figure}

The key element of the proof of  
is a tree-like recursive ordering over $\mbn^{d}$. To define it we introduce the following notation. For $t=(t_1,\ldots,t_d)\in \mbn^d$ define the "dyadic projection at scale $k$" as $\pi_k(t):=\paren{\lfloor 2^{-k}t_i\rfloor 2^k}_{1\leq i \leq d} \in 2^k\mbn^d$. Observe that $\pi_0(t)=t$, and that $\pi_k(t)=\bm{0}$ for $k\geq \log_2 \norm{t}_\infty$. 
Let $\leqlex$ denote the lexicographical order on $\mbn^d$. Denote $\lelex$ to be the associated strict order relation.
For two elements $t,t'$ of $\mbn^d$, define
\begin{equation} 
	\label{eq:defkappa_01}
	\kappa(t,t') = \min \{ k \in \mbn: \pi_k(t)=\pi_k(t')\}-1.
\end{equation}
Note that $\kappa(t,t')$ is always well-defined, since $\pi_k(t)=\pi_k(t')=\bm{0}$ for $k\geq \max(\norm{t}_\infty,\norm{t'}_{\infty})$, hence
the minimum in Equation~\eqref{eq:defkappa_01} is over a non-empty set. Furthermore $\kappa(t,t')=-1$ iff $t=t'$.  We define the following order $\preceq$ on $\mbn^d$
\begin{equation}
	\label{eq:deforder_01}
	t \preceq t' \;\;\; \text{ iff } \;\;\; \text{ either } \kappa(t,t')=-1 \text{ or } \pi_{\kappa(t,t')}(t) \lelex \pi_{\kappa(t,t')}(t').
\end{equation}
It is straightforward to check that $\preceq$ is a total order over $\mbn^d$ since $\leqlex$ is a total order over $\mbn^d$. This order can be described as the co-lexicographical order for the (one-to-one) sequence representation $(\pi_k(t))_{k\geq 0}$ of $t\in \mbn^d$, where the base order for the elements of the sequence is the usual lexicographical order. Equivalently, this is the co-lexicographical order on the (infinite) binary representation $((\lfloor t_i2^{-k} \rfloor \text{ mod } 2)_{i=d,\ldots,1})_{k\geq 0}$, where the vectorization is along (reverse) dimension first, then along scale. 
See Figure~\ref{pic:order} an illustration of the order $\preceq$.

For $t \in \mbn^d$, let
\begin{align} 
	\label{eq:defpikprec_01}
	\Pi_k^{\prec}(t) &:= \{ t' \in \mbn^d : \pi_k(t') \prec \pi_k(t) \},\\
	\Pi_k(t) &:= \{ t' \in \mbn^d : \pi_k(t') \preceq \pi_k(t) \}, \label{eq:defpik_01}
\end{align}
where $t\prec t'$ indicates that $t$ is strictly less than $t'$ for the order $\preceq$. Next we need the following Lemma which describes properties of the order $\preceq$.
\begin{lemma}
	\label{lem:proporder_01}
	Let ordering $\leq_{cw}$ be the partial order on $\mbn^d$ such that for $t=\paren{t_{1},\ldots,t_{d}} \in \mbn^d$, $t'=\paren[1]{t'_{1},\ldots,t'_{d}} \in \mbn^d$ we say that $t\leq_{cw} t'$ iff $t_{i} \leq t'_{i}$ for all $i \in \irg{d}$. The following statements hold true:
	\begin{itemize}
		\item[(0)] For any $k,\ell$ such that $k \leq \ell$ it holds $\pi_{\ell} \circ \pi_{k} = \pi_{k} \circ \pi_{\ell} = \pi_{\ell}$. 
		\item[(i)] The partial order $\leqcw$
		is compatible with both the total orders $\leqlex$ and $\preceq$, meaning that
		\[t \leqcw t' \;\; \Longrightarrow
		\;\; t \leqlex t' \text{ and }  t \preceq t'.
		\]
		\item[(ii)] For any $t\in\mbn^d$ and $k,\ell\in\mbn$ with $k \geq \ell$ it holds $\pi_k(t) \preceq \pi_\ell(t)$. In particular, in case $\ell=0$ it holds $\pi_k(t) \preceq t$.
		\item[(iii)] All applications $\pi_k$ are monotone nondrecreasing with respect to $\preceq$:
		\[
		\forall k \in \mbn, \;\; \forall t,t' \in \mbn^d: \qquad t \preceq t' \Rightarrow \pi_k(t) \preceq \pi_k(t').
		\]
		\item[(iv)] For an 
        integer $k$ put  $\cC_{k,0} := \irg{2^{k}}_0^d$, and, for $b\in 2^k\mbn^d$, put $\cC_{k,b} := \set{b} + \cC_{k,0}$.
		For any $t\in \cC_{k,b}$,
		it holds $\pi_k(t)=\pi_k(b)=b$,
		$\Pi^\prec_{k}(t) = \Pi^\prec_{k}(b)$ and $\Pi_{k}(t) = \Pi_{k}(b)$. 
		\item[(v)] For any $t\in\mbn^d$ and $k\in\mbn$, it holds
		\begin{align}
			\Pi_k^{\prec}(t) &= \{ t' \in \mbn^d : t' \prec \pi_k(t) \}, \label{eq:proppik1_01}\\
			\Pi_k(t) & = \Pi_k^{\prec}(t) \cup \cC_{k,\pi_k(t)}. \label{eq:proppik2_01}
		\end{align}
		\item[(vi)] For any $k \in \mbn_{>0}$ and $t \in \mbn^d$, 
		it holds
		\begin{equation}
			\label{eq:inclscales_01}
			\Pi_{k}^\prec(t) \subseteq \Pi^\prec_{k-1}(t) \subseteq \Pi_{k-1}(t) \subseteq \Pi_{k}(t).
		\end{equation}
	\end{itemize}
\end{lemma}

\begin{proof}
	$\mathbf{(0)}$ For $\ell \geq k$, the equality $\pi_{k} \circ \pi_{\ell} = \pi_{\ell}$ follows directly from the
	fact that $\lfloor 2^{\ell -k} \lfloor 2^{-\ell}t\rfloor \rfloor = 2^{\ell-k} \lfloor 2^{-\ell}t \rfloor$,
	while  the equality $\pi_{\ell} \circ \pi_{k} = \pi_{\ell}$ follows from $\lfloor 2^{k -\ell} \lfloor 2^{-k}t\rfloor \rfloor = \lfloor 2^{-\ell}t \rfloor$ which is easy to check.\\
	
	$ \mathbf{(i)}$ The implication for the lexicographical order is obvious;
	concerning the order $\preceq$, note that obviously all the mappings $\pi_k$ for $k\geq 0$ are non-decreasing for the
	partial order $\leqcw$, i.e. $t \leqcw t'$ implies $\pi_k(t) \leqcw \pi_k(t')$, in turn
	implying $\pi_k(t) \leqlex \pi_k(t')$ for all $k$, which finally entails $t \preceq t'$ from the definition. \\ 
	
	$\mathbf{(ii)}$ The claim follows directly from $\mathbf{(i)}$ since
	it follows from the definition that $\pi_k(t) \leqcw \pi_\ell(t)$
	if $k \geq \ell$.  \\
	
	$\mathbf{(iii)}$ Assume $t \prec t'$ and	let $\kappa = \kappa(t,t') \geq 0$. Then by definition
	of the order $\preceq$, for $\ell=\kappa+1$ it holds $\pi_\ell(t)=\pi_\ell(t')$, and therefore also further for any $\ell>\kappa$,
	since $\pi_\ell = \pi_\ell \circ \pi_{\kappa+1}$ by point $\mathbf{(0)}$.
	Thus, for $k>\kappa$
	it holds $\pi_k(t) \preceq \pi_k(t')$.
	On the other hand, it holds $\pi_\kappa(t) \lelex \pi_\kappa(t')$ and
	for $\ell<\kappa$, $\pi_\ell(t) \neq \pi_\ell(t')$.  
	For $k\leq \kappa$, put $u:=\pi_k(t), u':=\pi_{k}(t')$, then
	for any $\ell$ we have by point $\mathbf{(0)}$: $\pi_\ell(u) = \pi_{\max(\ell,k)}(t)$ and $\pi_\ell(u') = \pi_{\max(\ell,k)}(t')$.
	It follows that $\kappa(u,u')=\kappa$ and that
	the conditions for $u \prec u'$ are met.
	In both cases we have $\pi_k(t) \preceq \pi_k(t')$.\\
	
	$\mathbf{(iv)}$ For $u \in 2^k\mbn$, it holds $\lfloor 2^{-k}(u+v) \rfloor =u$
	iff $v \in \irg{2^k}_0$. It follows that for $t,t' \in \mbn^d$, $\pi_k(t) = \pi_k(t')$ iff $t' \in \cC_{k,\pi_k(t)}$. The claims follow from the definitions of $\pi_k,\Pi_k^{\prec}$ and $\Pi_k$.\\
	
	$\mathbf{(v)}$ For $t,t' \in \mbn^d$, if $t' \prec \pi_k(t)$ then $\pi_k(t') \preceq t' \prec \pi_k(t)$, from point $\mathbf{(ii)}$. Conversely, if $t' \succeq \pi_k(t)$, then
	$\pi_k(t') \succeq \pi_k(\pi_k(t)) = \pi_k(t)$, by $\mathbf{(iii)}$.
	Hence $t' \prec \pi_k(t)$ iff $\pi_k(t') \prec \pi_k(t)$. This establishes~\eqref{eq:proppik1_01}. Concerning~\eqref{eq:proppik2_01}, we have seen above
	(see proof of $\mathbf{(iv)}$) that  $\{ t' \in \mbn^d: \pi_k(t) = \pi_k(t')\}= \cC_{k,\pi_{k}\paren{t}}$, therefore
	\[
	\Pi_k(t)  = \Pi_k^{\prec}(t) \cup \{ t' \in \mbn^d: \pi_k(t) = \pi_k(t')\}  = \Pi_k^{\prec}(t) \cup \cC_{k,\pi_k(t)}.
	\]

	$\mathbf{(vi)}$ 
	It holds $\pi_{k-1}(t) \succeq \pi_k(t)$ from $\mathbf{(ii)}$.
	Then from~\eqref{eq:proppik1_01}, we deduce the inclusion $\Pi_{k}^\prec(t) \subseteq \Pi^\prec_{k-1}(t)$. The inclusion $\Pi^\prec_{k-1}(t) \subseteq \Pi_{k-1}(t)$
	is immediate from the definitions~\eqref{eq:defpikprec_01},~\eqref{eq:defpik_01}.
	Finally, for any $t' \in \Pi_{k-1}(t)$, by definition $\pi_{k-1}(t') \preceq \pi_{k-1}(t)$, so
	by $\mathbf{(iii)}$ and $\pi_k \circ \pi_{k-1} = \pi_k$, it holds $\pi_{k}(t') \preceq \pi_k(t)$, hence $t' \in \Pi_{k}(b)$, proving the last inclusion.  
\end{proof}

\begin{remark}  
	The choice of the lexicographical order in the definition~\eqref{eq:deforder_01} is largely arbitrary; any total order on $\mbn^d$ that is compatible with the coordinate-wise partial order would work, since it would result in the same properties as above, which are the only ones we will be using in the sequel. 
\end{remark}

For $t \in \mbn^d$ and an integer $k$, define
\begin{align}
	\begin{aligned}
		\label{eq:sig_algbr_01}
		\cF^{\prec}_{k}(t) &:= \fS(X_{t'}, t' \in \Pi_k^\prec(t)), \qquad & 
		\cF_{k}(t) 	&:= \fS(X_{t'}, t' \in \Pi_k(t)),
	\end{aligned}
\end{align}
where $\Pi_k^\prec,\Pi_k$ are as defined in~\eqref{eq:defpikprec_01},~\eqref{eq:defpik_01} (and $\fS(\emptyset)$ is the trivial $\sigma$-algebra).

For every element $t \in \mbn^d$, using the fact that $\Pi^{\prec}_k(t) = \emptyset$ for $k>\log_2 t$, we write the decomposition 
\[
X_{t}  - \e{X_{t}} =  \paren[2]{X_{t} - \e[1]{X_{t}|\cF^\prec_{0}\paren{t}}} + \sum_{k=1}^{\lfloor \log_2 \norm{t}_\infty \rfloor +1} \paren{\e[1]{X_{t}|\cF^\prec_{k}\paren{t}} -
	\e[1]{X_{t}|\cF^\prec_{k-1}\paren{t}}}.
\]

For any finite subset $A\subset \mbn^d$, denoting $\pi_k(A) = \{ \pi_k(t), t\in A\} \subset 2^k\mbn^d$ and $\norm{A}_\infty = \max_{t \in A} \norm{t}_\infty$, we have $A = \biguplus_{b \in \pi_k(A)} (A\cap \cC_{k,b})$, hence:
\begin{align}  
	S_{A}
	& = \sum_{t \in A} \paren{X_t-\e{X_t}} \notag \\
	&  = \sum_{t \in A} \paren[2]{X_{t} - \e[1]{X_{t}|\cF^\prec_{0}\paren{t}}}
	+ \sum_{k=1}^{\lfloor \log_2 \norm{A}_\infty \rfloor +1} \sum_{b\in \pi_{k}(A)} \sum_{t \in \cC_{k,b} \cap A} \paren{\e[1]{X_{t}|\cF^\prec_{k}\paren{t}} -
		\e[1]{X_{t}|\cF^\prec_{k-1}\paren{t}}} \notag \\
	& = \sum_{k=0}^{\lfloor \log_2 \norm{A}_\infty \rfloor+1} \sum_{b\in \pi_{k}(A)} Z_{b,k}(A),
	\label{eq:mart_dcmp_01}
\end{align}
where
\begin{align} \label{def:martz1_01}
	Z_{t,0}(A)
	& := {X_{t} - \e[1]{X_{t}|\cF^\prec_{0}\paren{t}}} ;\\
	\text{ and for } k \geq 1: \;\;\; Z_{b,k}(A)
	&:= \sum_{t \in \cC_{k,b} \cap A} \paren{\e[1]{X_{t}|\cF^\prec_{k}\paren{t}} -
		\e[1]{X_{t}|\cF^\prec_{k-1}\paren{t}}}. \label{def:martz2_01}
\end{align}
\begin{lemma}
	\label{lem:mart}
	Let $A\subset \mbn^d$ be a finite set. Let $k$ be a fixed integer. Then $(Z_{b,k}(A),\cF_{k}(b))_{b \in \pi_k (A)}$ is a martingale difference, where $\pi_k(A)$ is ordered by the total order $\preceq$ defined by~\eqref{eq:deforder_01}.
\end{lemma}

\begin{proof}
	We start with the special case $k=0$. In this case, since $\pi_0(t)=t$, we have $\cF^\prec_{0}\paren{t} = \fS(X_{t'}, t' \prec t)$, and
	$\cF_{0}\paren{t} = \fS(X_{t'}, t' \preceq t)$. It is straightforward that
	$Z_{0,t}(A)$ is $\cF_0(t)$-measurable, and that for any $t' \prec t$  we have $\cF_0(t') \subseteq \cF_0^\prec(t)$ thus
	$\e{Z_{0,t}(A)|\cF_0(t')}=0$; hence the claim. Let $k\geq 1$ be a fixed integer. The claim for $(Z_{b,k})$ relies on
	points $\mathbf{(iv)}$ and $\mathbf{(vi)}$ of Lemma~\ref{lem:proporder_01}, which straightforwardly imply for any $t\in \cC_{k,b}$ that  $\cF_{k}^\prec(b) = \cF_{k}^\prec(t) \subseteq \cF^\prec_{k-1}(t) \subseteq \cF_{k}(t) = \cF_{k}(b)$. Thus, $Z_{b,k}(A)$ is $\cF_k(b)$-measurable, and for any $b'\prec b$, since $\cF_{k}(b') \subseteq \cF^\prec_{k}(b)$,
	it holds $\e{Z_{b,k}(A)|\cF_{k}(b')}=0$, implying the claim.
\end{proof}\\


As announced at the beginning of the proof, the role of excluding elements
from the "frame" as constructed in Section~\ref{subsec:FDIS} is to
ensure that the remaining elements are sufficiently ``separated'' from the cell boundaries.
We recall that $\Lambda_{k,\delta} \bydef 2^{k} \mbn_{k>0} + \irg{\delta}_0 $. We need the following supporting result which estimates the distance  from any element of the set $\fF_{\bmd}$ to the boundaries of the cells containing it.

\begin{lemma}
	\label{lem:separ}
	For any  $k \in \mbn$, $\delta \in \irg{2^k}_0$ and $t\in(\mbn\setminus \Lambda_{k,\delta})^d$, it holds
	\[
	d_\infty(t, \Pi^\prec_k(t)) \geq \delta + 1. 
	\]
\end{lemma}
\begin{proof}
	Point $\mathbf{(i)}$ from Lemma~\ref{lem:proporder_01} implies that $\preceq$ is compatible with the partial coordinate-wise order $\leqcw$.
	This implies in particular that any $t'$ such that $\pi_k(t) \leqcw t'$ satisfies
	$\pi_k(\pi_k(t)) = \pi_k(t) \preceq \pi_k(t')$, and thus cannot belong to $\Pi^\prec_k(t)$. Therefore, for any $t' \in \Pi^\prec_k(t)$, there exists a coordinate $i$ such that $t'_i < \pi_k(t_i)$. In particular, $\pi_k(t_i)>0$, hence $\pi_k(t_i) \in 2^k\mbn_{>0}$. On the other hand, if we assume $t \in (\mbn\setminus \Lambda_{k,\delta})^d$ then $t_i \in \mbn\setminus (2^k\mbn_{>0} + \irg{\delta}_0)$. Since $t'_{i}<\pi_k(t_i) \leq t_i$, it must hold $t_i - t'_i \geq t_i - \pi_k(t_i) +1 \geq \delta+1$, implying the claim.
\end{proof}

\smallskip

We how have all ingredients to establish Proposition~\ref{prop:framedestimate_01}.

\begin{proof}[of~Proposition~\ref{prop:framedestimate_01}]
	We use the decomposition~\eqref{eq:mart_dcmp_01} with $A=\cR \cap \fF_{\delta}$, so that by the triangle inequality
	\begin{equation}
		\label{eq:dcpnormframed_01}
		\norm{S_{\cR \cap \fF_{\bmd}}}_p \leq \sum_{k=0}^{m(\cR)+1} \norm[4]{\sum_{b\in \pi_{k}(\cR \cap \fF_{\bmd})} Z_{b,k}(\cR \cap \fF_{\bmd})}_p,   
	\end{equation}
	where $Z_{b,k}(\cR \cap \fF_{\bmd})$ is defined in~\eqref{def:martz1_01},~\eqref{def:martz2_01}.
	
	We now estimate the norm of the martingale increments $Z_{b,k}(\cR \cap \fF_{\bmd})$ using Assumption~$\WD(p)$. We will denote below $Z_{b,k}=Z_{b,k}(\cR \cap \fF_{\bmd})$ and $S_k=\pi_k(\cR \cap \fF_{\bmd})$  to lighten notation.
	As a direct consequence of Lemma~\ref{lem:separ}, for any $t\in \fF_\bmd$ it holds $\cF_k^\prec(t) \subset \cM_{t,\delta_k+1}$ (as defined in Assumption~$\WD(p)$).
	Therefore, using this property and Jensen's inequality for $k=0$ we get:
	\begin{align}
		\norm{Z_{b,0}}_p
		& = \norm{X_{b} - \e[1]{X_{b}|\cF^\prec_{0}\paren{b}}}_p \notag \\
		& \leq \norm{\ee{}{X_{b}} - \ee{}{X_{b} | \cF^\prec_{0}\paren{b} }}_{p} + \norm{X_{b} - \ee{}{ X_{b}} }_{p} \notag \\
		&    \leq \norm{X_{b} - \e{X_b}}_p + \norm{\e[1]{X_{b}|\cM_{b,1}} - \e{X_b}}_p \notag\\
		& \leq    M_p (\varphi_{p}(0) + \varphi_{p}(1)),   \label{eq:term0_01}
	\end{align}
	while for $k\geq 1$:
	\begin{align}
		\norm{Z_{b,k}}_p
		& = \norm[3]{\sum_{t \in \cC_{k,b} \cap \cR \cap \fF_{\bmd}} \paren{\e[1]{X_{t}|\cF^\prec_{k}\paren{t}} -
				\e[1]{X_{t}|\cF^\prec_{k-1}\paren{t}}}}_p \notag \\
		& \leq \sum_{t \in \cC_{k,b} \cap \cR \cap \fF_{\bmd} } \paren{\norm[2]{ \e[1]{X_{t}|\cF^\prec_{k}\paren{t}} - \e{X_t}}_p + \norm[2]{
				\e[1]{X_{t}|\cF^\prec_{k-1}\paren{t}}-\e{X_t}}_p}\notag \\
		& \leq \sum_{t \in \cC_{k,b} \cap \cR \cap \fF_{\bmd} } \paren{\norm[2]{ \e[1]{X_{t}|\cM_{t,\delta_{k}+1}} - \e{X_t}}_p + \norm[2]{
				\e[1]{X_{t}|\cM_{t,\delta_{k-1}+1}}-\e{X_t}}_p}\notag \\
		& \leq \abs{\cC_{k,b} \cap \cR}  M_p \paren[1]{ \varphi_{p}(\delta_k+1) + \varphi_{p}(\delta_{k-1}+1)}.
		\label{eq:termk_01}
	\end{align}
	Note that we can subsume~\eqref{eq:term0_01} into~\eqref{eq:termk_01} by putting formally $\delta_{-1}:=-1$.
	By Lemma~\ref{lem:mart}, the sequence $(Z_{b,k},\cF_{k}(b))_{b \in S_k}$ is a martingale
	difference sequence (over $b$ for fixed $k$), therefore for $p\in[2,\infty)$ we apply
 the Marcinkiewicz-Zygmund inequality for martingales, with optimal constant obtained from Theorem 2.1 of \cite{Rio:09} 
 (see also Theorem 4.3 in \cite{Pinelis:94} for a related result under the assumption of $\cF_{k}(b)$-conditionally symmetric distributed martingale difference sequence) which, together with the nonincreasing character of $\varphi_{p}(\cdot)$
 {and nondecreasing character of $(\delta_k)_{k \geq 1}$} implies for $p\geq 2$:
 
	\begin{align}
		\norm[3]{\sum_{b\in S_k} Z_{b,k}}_p
		& \leq 2\sqrt{p} \paren[3]{ \sum_{b\in S_k}
			\norm{Z_{b,k}}_{p}^2 }^{\frac{1}{2}} \notag\\
		& \leq 2 \sqrt{p} M_p \paren[1]{ \varphi_{p}(\delta_k+1) + \varphi_{p}(\delta_{k-1}+1)}
		\paren[3]{ \sum_{b\in \pi_k(\cR)} \abs{\cC_{k,b} \cap \cR}^2 }^{\frac{1}{2}} \notag \\
		& \leq 4 \sqrt{p}M_{p}\varphi_{p}\paren{\delta_{k-1}+1}\paren[3]{ \sum_{b\in \pi_k(\cR)} \abs{\cC_{k,b} \cap \cR}^2 }^{\frac{1}{2}}.
		\label{eq:estnormz_01}
	\end{align}
	We now concentrate on the estimate for $\sum_{b\in \pi_k(\cR)} \abs{\cC_{k,b} \cap \cR}^2$. Put $q_i:= \big\lfloor \frac{N_i}{2^k} \big\rfloor$
	and $r_i := N_i- q_i2^k$, for $i \in \irg{d}$. Observe that $\pi_k(\cR) = \prod_{i=1}^d (2^k \irg{q_i+1}_0)$;
	for $b=(b_1,\ldots,b_d) \in \pi_k(\cR)$, the set $\cC_{k,b} \cap \cR$ is a hyperrectangle with side-lengths:
	\[
	\ell_k(b_i,N_i) :=
	\begin{cases}
		2^k & \text{ if } b_i < 2^kq_i;\\
		r_i & \text{ if } b_i = 2^kq_i.
	\end{cases}
	\]
	Hence, it holds
	\begin{align}
		\sum_{b\in \pi_k(\cR)} \abs{\cC_{k,b} \cap \cR}^2
		= \sum_{ b \in \prod_{i=1}^d (2^k \irg{q_i+1}_0)} \prod_{i=1}^d \ell_k(b_i,N_i)^2
		& = \prod_{i=1}^d \sum_{j=0}^{q_i} \ell_k(2^kj,N_i)^2 \notag \\
		& = \prod_{i=1}^d (q_i2^{2k} + r_i^2) \notag \\
		& \leq \prod_{i=1}^d \paren[1]{(N_i-r_i)\min(2^{k},N_i) + r_i \min(2^k,N_i)} \notag\\
		& = \prod_{i=1}^d \paren[1]{N_i\min(2^{k},N_i)} \notag \\
		& = \abs{R} \abs{R \cap \cC_{0,k}}. \label{eq:estsumsqcard_01}
	\end{align}
	The claimed estimate for $2\leq p< \infty$ follows by using~\eqref{eq:estnormz_01} and~\eqref{eq:estsumsqcard_01} into~\eqref{eq:dcpnormframed_01} and straightforward computations.
	In the case of $p=\infty$, we can apply the bounded martingale difference inequality
	\cite{Azuma:67} stating that the sum 
	$\sum_{b \in S_{k}} Z_{b,k}$ is sub-Gaussian such that $ \norm{\sum_{b \in S_{k}}Z_{b,k}}_{\SG}  \leq \paren[1]{\norm[1]{\sum_{b\in S_k} Z_{b,k}^2}_\infty}^{\frac{1}{2}}$ and using the triangle inequality for the sub-Gaussian norm 
 \hfill \break $\norm[1]{\sum_{k=0}^{m\paren{\cR} +1} {\sum_{b \in \pi_{k}\paren{\cR \cap \fF_{\bmd}} }}Z_{b,k}\paren{\cR \cap \fF_{\bmd}} }_{\SG}$ over scales $k \in \irg{m\paren{\cR}+2}_{0}$. 
	All other arguments are as in the case $p<\infty$. 
\end{proof}

	{\bf Acknowledgements.} The authors want to thank J\'erome Dedecker for very insightful discussions and comments and
 Andrey Pilipenko for equally interesting exchanges. The work of the three authors was partially supported by the DFG CRC 1294 'Data Assimilation', Project A03. The work of A. Carpentier is partially supported by the Deutsche Forschungsgemeinschaft (DFG) Emmy Noether grant MuSyAD (CA 1488/1-1), by the DFG Forschungsgruppe FOR 5381 "Mathematical Statistics in the Information Age - Statistical Efficiency and Computational Tractability", Project TP 02, by the Agence Nationale de la Recherche (ANR) and the DFG on the French-German PRCI ANR ASCAI CA 1488/4-1 "Aktive und Batch-Segmentierung, Clustering und Seriation: Grundlagen der KI". GB acknowledges support of the ANR under ANR-19-CHIA-0021-01 ``BiSCottE'', and IDEX REC-2019-044. Oleksandr Zadorozhnyi acknowledges Alexander von Humboldt Foundation (Research Group Linkage cooperation Singular diffusions: analytic and stochastic approaches between the University of Potsdam and the Institute of Mathematics of the National Academy of Sciences of Ukraine).


\bibliography{references}
\bibliographystyle{imsart-nameyear.bst}
\end{document}